\documentclass[a4paper,12pt]{article}
\usepackage{amsmath,amsfonts,amssymb,bm,amsthm}
\usepackage{mathrsfs}
\usepackage{cases}
\usepackage{cite}
\usepackage{indentfirst}
\usepackage[usenames]{color}
\usepackage{hyperref}
\usepackage{fancyhdr}
\usepackage{geometry}
\definecolor{webgreen}{rgb}{0,.5,0}
\definecolor{webbrown}{rgb}{.6,0,0}

\usepackage{color}

\newtheorem{theorem}{Theorem}

\newtheorem{definition}{Definition}

\newtheorem{corollary}{Corollary}

\allowdisplaybreaks

\begin{document}

\title{Analogy between geodesic equation and the GCHS on Riemannian manifolds}
\vspace{.3 cm}
\author{Gen Wang}
\date{{\em \small School of Mathematical Sciences, Xiamen University,\\
     Xiamen, 361005, P.R.China.}}

\maketitle
\abstract{Enlightened by the similar equation form between the GCHS \footnote{GCHS: Generalized Covariant Hamilton System\\GSPB:Generalized structural Poisson bracket} defined by the GSPB and the geodesic equation expressed by geospin variable,  we find a deep connection between the geospin matrix and S-dynamics.
In this contrastive way, we actually prove that the GCHS is a compatible theory suitable for the curved spacetime as primitively stated.  By contrast, geospin matrix in Riemannian geometry has the same physical nature as S-dynamics in GCHS. We obtain a fact that geodesic equation can be naturally derived by the GCHS in terms of the velocity field.

We strictly prove that the geometrio
$\hat{S}{{\left( {{x}_{k}},{{p}_{i}},H \right)}^{T}}={{\left( {{b}_{k}},{{A}_{i}},w \right)}^{T}}$
holds by using structural operator $\hat{S}$ directly induced by structural derivative ${A}_{i}$ in terms of position ${x}_{k}$, momentum ${p}_{i}$ and Hamiltonian $H$ respectively. It evidently proves that the GCHS on the Riemannian manifold is certainly determined by the Christoffel symbols.  As an application, we consider the GCHS on Riemannian geometry.

}

\newpage
\tableofcontents

\section{Introduction}
\label{sec:intro}

\subsection{Overviews of the GCHS and geodesic equation}
In this part, we firstly give a simple overviews of the GCHS \cite{1} and geodesic equation \cite{2,3,4,5,6,7,8,9} for the convenience of later discussion.

We begin with the generalized Poisson bracket (GPB) that is defined as the bilinear operation
$${{\left\{ f,g \right\}}_{GPB}}={{J}_{ij}}\frac{\partial f}{\partial {{x}_{i}}}\frac{\partial g}{\partial {{x}_{j}}}$$ where structural matrix $J$ satisfies antisymmetric ${{J}_{ij}}={{\left\{ {{x}_{i}},{{x}_{j}} \right\}}_{GPB}}=-{{J}_{ji}}$.
Generalized Hamiltonian system (GHS) is defined by generalized Poisson bracket (GPB)
\begin{equation}\label{eq4}
  {\dot{x_{i}}}=\frac{dx_{i}}{dt}=\left\{ {{x}_{i}},H \right\}_{GPB}={{{J}_{ij}}\frac{\partial H}{\partial {{x}_{j}}}},~~~x\in {{\mathbb{R}}^{m}}
\end{equation}
Subsequently, \cite{1} has further developed generalized structural Poisson bracket (GSPB) equipped with geometric bracket $G\left( s ;f,H \right)$ (in this case for functions $f,H$) to define generalized covariant Hamilton system (GCHS) as a complete version of generalized Hamiltonian system, which is more concise and complete, thus promoting the further development of generalized Hamilton system dynamics,  its basic theoretical framework is perfect, and it remains to be invariant.  Above all,
\[GPB\to GSPB,~~~~GHS\to GCHS\]
 The GCHS is precisely given by
\[\frac{\mathcal{D}f}{dt}=\left\{ f,H \right\}={{\left\{ f,H \right\}}_{GPB}}+G\left( s ;f,H \right)=\frac{df}{dt}+wf\]
for all $f\in {{C}^{\infty }}\left( M,\mathbb{R} \right)$, where $s$ is a real structure function on $M$. Accordingly,  the covariant equilibrium equation of the GCHS is given by $\frac{\mathcal{D}f}{dt}=\left\{ f,H \right\}=0$ which means that $f$ is a covariant conserved quantity.

The geodesic equation is a second-order ODE which correspondingly has a unique solution, given an initial position and an initial velocity. Indeed, the equation means that the acceleration vector of the curve has no components in the direction of the surface \cite{2,4,6}. So, the motion is completely determined by the bending of the surface. This is also the idea of general relativity where particles move on geodesics and the bending is caused by the gravity. It also can be written in the form \cite{3}
\begin{equation}
  \left\{ \begin{matrix}
   \frac{d{{x}^{k}}}{dt}={{v}^{k}}\begin{matrix}
   {} &  \\
\end{matrix}  \notag\\
   \frac{d{{v}^{k}}}{dt}+W_{j}^{k}{{v}^{j}}=0 \notag\\
\end{matrix} \right.\notag
\end{equation}
where $W_{j}^{k}$ are geospin variables.
The quantity on the left hand side of this equation is the acceleration of a particle, so this equation is analogous to Newton's laws of motion, which likewise provide formulae for the acceleration of a particle. In \cite{5}, we have proven that the geospin variable as a dynamical system is key element to form the geometrodynamics Riemannian manifolds based on the Cartan structural
equation.

The aim of this paper is to consider the following equations between the covariant equilibrium equation of the GCHS and geodesic equation together
\[\left\{ \begin{matrix}
   \frac{\mathcal{D}f}{dt}=\left\{ f,H \right\}={{\left\{ f,H \right\}}_{GPB}}+G\left( s ;f,H \right)=\frac{df}{dt}+wf=0  \\
   \frac{{{d}^{2}}{{x}^{k}}}{d{{t}^{2}}}+\Gamma _{ij}^{k}\frac{d{{x}^{j}}}{dt}\frac{d{{x}^{i}}}{dt}=0,~~i,j,k=1,\cdots ,n  \\
\end{matrix} \right.\]
It mathematically proves and states that covariant equilibrium equation of the GCHS which is formally similar to the geodesic equation, and indeed, a covariant form just as geodesic equation stated. They both are complete covariant form in a same mathematic pattern
\[\begin{matrix}
   \frac{\mathcal{D}f}{dt}=\left\{ f,H \right\}=\frac{df}{dt}+wf=0  \\
   \frac{dv}{dt}+Wv=0  \\
\end{matrix}\]where $W$ is the geospin matrix \cite{3}.
This analogy strongly implies that the GCHS is suitable for the general case, and it's covariant and complete one,  namely, it can apply to the Riemannian manifold, the non-Euclidean case. Since the same pattern can be built up by this analogical method. This provides an interpretation of what the S-dynamics $w$ and the geospin matrix $W$ are, they behave like a similar behavior.

This paper is organized as follows. The section 1 simply reviews the GCHS and geodesic equation, how to analogously study them together, and the geospin matrix and geospin variables introduced by \cite{3} are given for geodesic equation. In section 2, we will prove a fundamental result only concerned with the structural operator induced by structural derivative.
The section 3 gives some results of the GCHS on the Riemannian manifolds, including the geomrtrio and acceleration.  In section 4, we take stock of the closely connection between the GCHS and geodesic equation, then it actually proves that the GCHS is a compatible theory suitable for the curved spacetime as \cite{1} primitively stated, and it concludes that geospin variables and S-dynamics describe the same mathematical, physical phenomena.

\subsection{GSPB, GCHS, geometric bracket}
In order to better understand the concept of the GCHS for the subsequent discussions, in this section, we briefly retrospect some basic concepts of entire framework of GCHS \cite{1}.  Let $M$ be a smooth manifold and let $s$ be a smooth real structural function on $M$ which is completely determined by the structure of manifold $M$.
\begin{theorem}[GSPB]\label{t1}\cite{1}
The GSPB on $M$ of two functions $f,g\in {{C}^{\infty }}\left( M,\mathbb{R} \right)$ is defined as
\[\left\{ f,g \right\}={{\left\{ f,g \right\}}_{GPB}}+G\left( s ;f,g \right) \]for ${{\mathbb{R}}^{r}}$, the geometric bracket is
\begin{align}
 G\left( s,f,g \right) =f{{\left\{ s ,g \right\}}_{GPB}}-g{{\left\{ s ,f \right\}}_{GPB}} \notag
\end{align} in which $s$ is structural function associated with $M$.
\end{theorem}
Intuitively, one can understand this defined formulation by noting that the second part, geobracket, can be partially treated as the second part of the geodesic equation, and in so doing will make it clear. The GSPB depends smoothly on both ${{\left\{ f,g \right\}}_{GPB}}$ and the geobracket $G\left( s ;f,g \right) $. Obviously, the geobracket $G\left( s ;f,g \right) $ is necessary for a complete Hamiltonian syatem, as a result of the geobracket $G\left( s ;f,g \right) $, it can be generally used to depict nonlinear system in a non-Euclidean space, and for general manifolds.   As GSPB defined, the GCHS is completely assured and determined by the GSPB. Therefore, the GCHS can be distinctly obtained as follows.

\begin{theorem}[GCHS]\cite{1}
The GCHS on $M$ of two functions $f,H\in {{C}^{\infty }}\left( M,\mathbb{R} \right)$ is defined as
\[ \frac{\mathcal{D}f}{dt}=\left\{ f,H \right\}={{\left\{ f,H \right\}}_{GPB}}+G\left( s ;f,H \right) \]for ${{\mathbb{R}}^{r}}$,  where geobracket in terms of $f$ and Hamiltonian $H$ is
\begin{align}
 G\left( s ;f,H \right) =f{{\left\{ s ,H \right\}}_{GPB}}-H{{\left\{ s ,f \right\}}_{GPB}} \notag
\end{align}
\end{theorem}In fact, general covariance of the GCHS holds naturally.
The GCHS is totally dependent on the GSPB and accordingly defined. Notice that the GCHS is covariant, it means that the GCHS is complete, and its form is invariant under the coordinate transformation, just like the covariance of the geodesic equation. So, the motion is completely determined by the structure of the manifold $M$. This is also like the idea of general relativity where particles move on geodesics and the bending is caused by the gravity.

\begin{theorem}[TGHS,S-dynamics, GCHS]\cite{1}
\begin{description}
\item[TGHS:] $\frac{df}{dt}= {\dot {f}}={{\left\{ f,H \right\}}_{GPB}}-H{{\left\{ s ,f \right\}}_{GPB}}$.
  \item[S-dynamics:] $\frac{ds }{dt}=w={{\left\{ s ,H \right\}}_{GPB}}=\left\{ 1,H \right\}$.
  \item[GCHS:] $\frac{\mathcal{D}f}{dt}=\left\{ f,H \right\}={{\left\{ f,H \right\}}_{GPB}}+G\left( s ;f,H \right)$.
\end{description}
\end{theorem}
By defining the GSPB which can basically solve and describe an entire Hamiltonian system, even nonlinear system. Thus, in order to state this point. To start with, we recall the basic notions of the GSPB. We explain that the GCHS corresponds to non-Euclidean space which corresponds to curved coordinate systems in curved space-time. The geobracket $G\left( s ;f,H \right)$ is the part that has been neglected for a long time which is rightly the correction term for the general covariance which certainly means nonlinear system.  It is precisely because of the appearance of the second term  $G\left( s ;f,H \right)$  that the GCHS can describe and explain the generalized nonlinear system, the nonlinear Hamiltonian system.

In the following, the covariant equilibrium equation of the GCHS as a special case is naturally generated as follows.
\begin{corollary}\cite{1}\label{c4}
 The covariant equilibrium equation of the GCHS is $\frac{\mathcal{D}f}{dt}=\left\{ f,H \right\}=0$ such that $${{\left\{ f,H \right\}}_{GPB}}+G\left( s ;f,H \right)=0$$holds, where $s$ is a structure function of the manifolds $M$, then $f$ is called covariant conserved quantity.
\end{corollary}
Obviously, the covariant equilibrium equation of the GCHS leads to a formal solution that is given by \[f\left( t \right)={{f}_{0}}{{e}^{-wt}}\] with initial condition ${{f}_{0}}=f\left( 0 \right)$, where we suppose that ${{\partial }_{t}}w=0$. More precisely, if $w>0$ is given, then there exists a bound \[f\left( t \right)={{f}_{0}}{{e}^{-wt}}\le {{f}_{0}}\]
Taking a Taylor expansion leads to \[f\left( t \right)={{f}_{0}}{{e}^{-wt}}={{f}_{0}}\left( 1-wt+\cdots  \right)={{f}_{0}}+{{f}_{0}}Q\left( w,t \right)\]where $Q\left( w,t \right)={{e}^{-wt}}-1$ is the correction term.

The critical points of the GCHS are precisely corresponding to the geodesics-like.
In an appropriate sense, zeros of the GCHS are thus regarded as function $f$ through geodesics-like.

\begin{corollary}\label{c3}
The generalized force field $F$ on the generalized Poisson manifold $\left( P,S,\left\{ \cdot ,\cdot  \right\} \right)$ is shown as  $$F=\frac{d}{dt}p=-DH$$ Its components is $F_{k}={\dot{p}_{k}}=-{{D}_{k}}H$.
\end{corollary}
Hence the GCHS can be naturally rewritten in the form
$$\frac{\mathcal{D}f}{dt}=\left\{ f,H \right\}={{J}_{ij}}{{D}_{i}}f{{D}_{j}}H={{J}_{ji}}{F_{j}}{{D}_{i}}f$$
 Based on the GSPB, theoretically, we have the conservation of energy given by
$\left\{ H,H \right\}={{J}_{ij}}{{F}_{i}}{{F}_{j}}=0$  in which ${F}_{i}=-{{D}_{i}}H$ is the components of generalized force $F$, and the TCHS is rewritten as $${\dot{x}_{k}}={{J}_{kj}}{{D}_{j}}H={{J}_{jk}}{{F}_{j}}$$
As concerned, by employing the generalized force field $F$, the GCHS also has the form given by
 \begin{align}
\frac{\mathcal{D}f}{dt} &=\left\{ f,H \right\}=-{{D}^{T}}fJF\notag \\
 & =-{{\nabla }^{T}}fJF-f{{\nabla }^{T}}sJF \notag\\
 & ={{J}_{ji}}{{D}_{i}}f{{F}_{j}} \notag\\
 & ={{J}_{ji}}{{F}_{j}}{{\partial }_{i}}f+f{{J}_{ji}}{{F}_{j}}{{\partial }_{i}}s \notag
\end{align}
Thusly, the TGHS and S-dynamics respectively are shown as
\begin{align}
  & df/dt=-{{\nabla }^{T}}fJF={{J}_{ji}}{{F}_{j}}{{\partial }_{i}}f \notag\\
 & ds/dt=w=-{{\nabla }^{T}}sJF={{J}_{ji}}{{F}_{j}}{{\partial }_{i}}s\notag
\end{align}
Hence, it leads to the ordinary time derivative
\begin{equation}\label{eq6}
  d/dt={{J}_{ji}}{{F}_{j}}{{\partial }_{i}}
\end{equation}

\begin{theorem}
There exists two different canonical Hamilton equations on manifolds respectively given by
  \begin{description}
    \item[Canonical thorough Hamilton equations] \[\frac{d{{x}_{k}}}{dt}={\dot{x}_{k}}={{J}_{kj}}{{D}_{j}}H,~~\frac{d{{p}_{k}}}{dt}={\dot{p}_{k}}=-{{D}_{k}}H\]

    \item[canonical covariant Hamilton equations]
  \[\frac{\mathcal{D}{{x}_{k}}}{dt}=\left\{{{x}_{k}} ,H \right\},~~\frac{\mathcal{D}{{p}_{k}}}{dt}=\left\{{{p}_{k}} ,H\right\}\]
  \end{description}
where $\left\{ \cdot ,~\cdot  \right\}={{\left\{ \cdot ,~\cdot  \right\}}_{GPB}}+G\left(s,\cdot ,~\cdot  \right)$ is generalized structure Poisson bracket, ${{D}_{i}}={{\partial }_{i}}+{{\partial }_{i}}s $ is generalized derivative operators.
\end{theorem}

\begin{corollary}\label{c2}
  By corollary \ref{c3}, canonical thorough Hamilton equations (CTHE) can be simply rewritten as
  \[\frac{d{{x}_{k}}}{dt}={\dot{x}_{k}}={{J}_{jk}}{F_{j}},
        ~~\frac{d{{p}_{k}}}{dt}={\dot{p}_{k}}=F_{k}\]
\end{corollary}
With the corollaries \ref{c3} and \ref{c2}, we have some derivations given by
\begin{corollary}
The ordinary time derivative \eqref{eq6} can be rewritten in a form
\[d/dt={\dot{x}_{i}}{{\partial }_{i}}\]$$\mathcal{D}/dt={\dot{x}_{i}}{{D }_{i}}$$
where ${\dot{x}_{i}}$ is the TGHS,
and then for the generalized force field $F$, we get
\[d{{p}_{k}}/dt={{F}_{k}}={\dot{x}_{i}}{{\partial }_{i}}{{p}_{k}}={\dot{x}_{i}}{{\xi }_{ik}}\]in terms of the momentum.
Similarly, for the S-dynamics, it yields
\begin{equation}\label{eq7}
  ds/dt=w={\dot{x}_{i}}{{\partial }_{i}}s={\dot{x}_{i}}{{A}_{i}}
\end{equation}in terms of the structure function $s$,
and the TGHS is given by
\[df/dt={\dot{x}_{i}}{{\partial }_{i}}f\]for any function $f$.
\end{corollary}
Note that the S-dynamics expressed by \eqref{eq7} implies the relation between the structure derivative ${{A}_{i}}$ and the S-dynamics itself, the S-dynamics is only induced by the structure derivative ${{A}_{i}}$, and based on the manifolds itself.

\subsection{Geodesic equation on Riemannian geometry}
Riemannian geometry \cite{3,7,8,9} is concerned with the objective entity which has nothing to do with all coordinates.   A smooth Riemannian manifold $(M,g)$ is a real smooth manifold $M$ equipped with an inner product $g_{p}$ on the tangent space $T_{p}M$ at each point $p$ that varies smoothly from point to point in the sense that if $X$ and $Y$ are differentiable vector fields on $M$, then $p\mapsto g(X(p),Y(p))$ is a smooth function. The family $g_{p}$ of inner products is called a Riemannian metric tensor.

\begin{theorem}\cite{3,7}
In a local coordinate system $\left( U,{{x}^{i}} \right)$, connection $\nabla\colon\mathfrak{X}\left( U \right)\times \mathfrak{X}\left( U \right)\to \mathfrak{X}\left( U \right)$ can be arbitrarily determined, law of motion at the natural shelf field $\left\{ \frac{\partial }{\partial {{x}^{i}}};1\le i\le m \right\}$:
\begin{equation}\label{eq5}
  {{\nabla}_{\frac{\partial }{\partial {{x}^{k}}}}}\frac{\partial }{\partial {{x}^{j}}}=\Gamma _{jk}^{i}\frac{\partial }{\partial {{x}^{i}}}
\end{equation}
\end{theorem}

\begin{definition}\cite{7,8,9}
 A smooth curve $\gamma =\left[ a,b \right]\to M$, which satisfies (with  ${\dot{x}^{i}}=\frac{d}{dt}{{x}^{i}}\left( \gamma \left( t \right) \right)$ etc.)
 \[\ddot{x}^{i}(t)+\Gamma_{j k}^{i}(x(t)) \dot{x}^{j}(t) \dot{x}^{k}(t)=0, \text { for } i=1, \ldots, n\]
 is called a geodesic equation.
\end{definition}

\subsection{Geospin variables for geodesic equation}
In this section, we simply review geospin variables for geodesic equation given by \cite{3},  it shows that the equation of geodesics can be simplified by geospin matrix, in the case of Riemannian and pseudo-Riemannian manifolds.

\begin{definition}\cite{3}
Let $(M,g)$ be a Riemannian manifold with Eq\eqref{eq5}, two kinds of geospin variables can be defined as
  \begin{equation}
 W_{i}^{j}=\Gamma _{ik}^{j}{{v}^{k}},~~~{{W}_{ik}}=\Gamma _{ik}^{j}{{v}_{j}}
\end{equation}
where  ${{v}^{j}}, {{v}_{j}}$ are the components of  $v={{v}^{j}}\frac{\partial }{\partial {{x}^{j}}}={{v}_{j}}\frac{\partial }{\partial {{x}_{j}}}=\frac{dx}{dt}$ respectively.
\end{definition}
As geospin variables defined above, we mainly study the (1,1) form geospin variable $W_{i}^{j}$ that relates to the geospin matrix $W$ below.

\begin{definition}\cite{3}
The geospin matrix can be defined as
  $$W=\left( W_{j}^{k} \right)$$
where $W_{j}^{k}$ are geospin variables.
\end{definition}
Obviously, the symmetric holds ${{W}_{kj}}={{W}_{jk}}$. The covariant derivative of a vector field $v$ can be rewritten as follows \cite{3}
\begin{equation}\label{eq11}
  {{\nabla }_{k}}{{v}^{j}}=\frac{\partial {{v}^{j}}}{\partial {{x}^{k}}}+W_{k}^{j},~~{{\nabla }_{k}}{{v}_{j}}=\frac{\partial {{v}_{j}}}{\partial {{x}^{k}}}-{{W}_{kj}}
\end{equation}
 As a consequence, the geodesic for Riemannian manifolds can be rewritten in a compact and simple form by using the geospin variables as follows:
\begin{equation}\label{eq3}
  \left\{ \begin{matrix}
   \frac{d{{x}^{k}}}{dt}={{v}^{k}}\begin{matrix}
   {} &  \\
\end{matrix}  \\
   \frac{d{{v}^{k}}}{dt}+W_{j}^{k}{{v}^{j}}=0 \\
\end{matrix} \right.
\end{equation} More precisely, geodesic equation Eq\eqref{eq3} based on new matrix variable is written in a matrix form as
\begin{equation}
  \left\{ \begin{matrix}
   \frac{d{{x}}}{dt}={{v}}\begin{matrix}
   {} &  \\
\end{matrix}  \notag \\
   \frac{d{{v}}}{dt}=-Wv   \notag\\
\end{matrix} \right.
\end{equation}
This is a dynamic system. This formulation of the geodesic equation of motion can be useful for computer calculations and to compare general relativity with Newtonian gravity. The geodesics equation is transformed to initial problem \cite{3}
\begin{equation}\label{eq14}
  \frac{dv}{dt}=-Wv,~~~v\left( {{t}_{0}} \right)={{v}_{0}}
\end{equation}
 The left side of the equation is the acceleration vector of the curve on the manifold, so the equation means that the geodesic line is a curve with zero acceleration on the manifold, so the geodesic line must be a constant velocity curve.

\section{Structural operator and geometrio}
In this section, by defining the structural operator induced by the structure function $s$ in GSPB, we get an important theorem related to the basic building block of the GCHS.

\begin{definition}[Structural operator]\label{d1}\cite{1}
The structural operator on manifold $M$ is defined as
\[\hat{S}\equiv {{A}^{T}}JD ={{J}_{ij}}{{A}_{i}}{{D}_{j}}={{J}_{ij}}{{A}_{i}}{{\partial }_{j}}=X_{s}\]
where $D=\nabla +A$ is generalized gradient operator, $X_{s}$ is structural vector field.
\end{definition}
Obviously, we have \[\widehat{S}s={{J}_{ij}}{{A}_{i}}{{D}_{j}}s={{X}_{s}}s={{J}_{ji}}{{A}_{i}}{{A}_{j}}=0\]
\[\widehat{S}\left( Const \right)={{J}_{ij}}{{A}_{i}}{{D}_{j}}\left( Const \right)={{X}_{s}}\left( Const \right)=0\]
With the definition \ref{d1}, we can nicely prove the geometrio as follows:
\begin{theorem}[Geometrio]\label{t2}
The structural operator $\widehat{S}$ can induce the following geometrio
 \[\hat{S}{{\left( {{x}_{k}},{{p}_{i}},H \right)}^{T}}={{\left( {{b}_{k}},{{A}_{i}},w \right)}^{T}}\]
in terms of position ${x}_{k}$, momentum ${p}_{i}$ and Hamiltonian $H$ respectively.
\begin{proof}
  Based on the structural operator $\hat{S}={{J}_{ij}}{{\partial }_{i}}s{{\partial }_{j}}$ in \ref{d1}, we get
\begin{align}
{{b}_{k}}  & ={{\left\{ s,{{x}_{k}} \right\}}_{GPB}}=\hat{S}{{x}_{k}} \notag\\
 & ={{J}_{ij}}{{\partial }_{i}}s{{\partial }_{j}}{{x}_{k}} \notag\\
 & ={{J}_{ij}}{{\partial }_{i}}s{{\delta }_{jk}} \notag\\
 & ={{J}_{ik}}{{\partial }_{i}}s \notag\\
 & ={{J}_{ik}}{{A}_{i}} \notag
\end{align}in terms of position ${x}_{k}$, thusly, the structural operator is rewritten as $\hat{S}={{b}_{j}}{{\partial }_{j}}$,
and
\begin{align}
 {{A}_{k}} & =\hat{S}{{p}_{k}}={{\left\{ s,{{p}_{k}} \right\}}_{GPB}}=-{{\left\{ {{p}_{k}},s \right\}}_{GPB}} \notag\\
 & ={{J}_{ij}}{{\partial }_{i}}s{{\partial }_{j}}{{p}_{k}}\notag \\
 & =-{{J}_{ji}}{{\partial }_{j}}{{p}_{k}}{{\partial }_{i}}s \notag\\
 & ={{J}_{ij}}{{\partial }_{j}}{{p}_{k}}{{\partial }_{i}}s \notag\\
 & ={{\delta }_{ik}}{{\partial }_{i}}s \notag\\
 & ={{\partial }_{k}}s \notag
\end{align}in terms of momentum ${p}_{k}$, and
\begin{align}
  w& =\hat{S}H={{\left\{ s,H \right\}}_{GPB}} \notag\\
 & ={{J}_{ij}}{{\partial }_{i}}s{{\partial }_{j}}H \notag\\
 & ={{J}_{ij}}{{A}_{i}}{{\partial }_{j}}H \notag\\
 & ={{J}_{ij}}\hat{S}{{p}_{i}}{{\partial }_{j}}H \notag\\
 & ={{b}_{j}}{{\partial }_{j}}H \notag
\end{align}in terms of Hamiltonian $H$.
Most important, their relation is \[{{b}_{k}}=\hat{S}{{x}_{k}} ={{J}_{ik}}\hat{S}{{p}_{i}}={{J}_{ik}}{{A}_{i}}\]
In conclusively, we can rewrite geometrio as form
\begin{align}\label{eq1}
 Geo\left( {{x}_{k}},{{p}_{i}},H \right) &=\hat{S}{{\left( {{x}_{k}},{{p}_{i}},H \right)}^{T}}\\
 & =\hat{S}\left( \begin{matrix}
   {{x}_{k}}  \\
   {{p}_{i}}  \\
   H  \\
\end{matrix} \right)=\left( \begin{matrix}
   {{b}_{k}}  \\
   {{A}_{i}}  \\
   w  \\
\end{matrix} \right) \notag\\
 & =\left( \begin{matrix}
   {{J}_{ik}}{{A}_{i}}  \\
   {{A}_{i}}  \\
   {{J}_{ij}}{{A}_{i}}{{\partial }_{j}}H  \\
\end{matrix} \right) =\left( \begin{matrix}
   {{J}_{ik}}  \\
   1  \\
   {{J}_{ij}}{{\partial }_{j}}H  \\
\end{matrix} \right){{A}_{i}}
\notag
\end{align}
where $T$ means transposition of the matrix.
\end{proof}
\end{theorem}
Obviously, we can see based on the \eqref{eq1} that geometrio is only induced by the structure derivative ${{A}_{i}}$ to form a complete theory that is GCHS, since geometrio represents how the manifolds work.  In other words,
\begin{align}
  Geo\left( {{x}_{k}},{{p}_{i}},H \right)&=\hat{S}{{\left( {{x}_{k}},{{p}_{i}},H \right)}^{T}}={{\left( {{b}_{k}},{{A}_{i}},w \right)}^{T}} \notag\\
 & ={{\left( {{J}_{ik}},1,{{J}_{ij}}{{\partial }_{j}}H \right)}^{T}}{{\partial }_{i}}s \notag
\end{align}
where structure function $s$ is only associated with the manifolds or the spacetime.

In fact,  by employing the generalized force field given by corollary \ref{c3}, the S-dynamics can be rewritten as
\begin{align}
  w& =\widehat{S}H={{J}_{ij}}{{A}_{i}}{{D}_{j}}H={{X}_{s}}H \notag\\
 & =-{{J}_{ij}}{{A}_{i}}{{F}_{j}}={{J}_{ji}}{{A}_{i}}{{F}_{j}} \notag
\end{align}
Furthermore,
\begin{align}
 \widehat{S}w &={{\widehat{S}}^{2}}H={{J}_{ij}}{{A}_{i}}{{D}_{j}}w={{X}_{s}}w \notag\\
 & ={{J}_{ij}}{{A}_{i}}{{\partial }_{j}}w={{b}_{j}}{{\partial }_{j}}w  \notag
\end{align}
And similarly, for the structure derivative ${{A}_{k}}$, it has
\begin{align}
 \widehat{S}{{A}_{k}} & ={{\widehat{S}}^{2}}{{p}_{k}}={{J}_{ij}}{{A}_{i}}{{D}_{j}}{{A}_{k}}={{X}_{s}}{{A}_{k}} \notag\\
 & ={{J}_{ij}}{{A}_{i}}{{\partial }_{j}}{{A}_{k}}={{b}_{j}}{{\partial }_{j}}{{\partial }_{k}}s  \notag
\end{align}

With the help of the theorem \ref{t2}, we can rewrite the geometric canonical commutation relation as follows.
\begin{corollary}
  The geometric canonical commutation relation is rewritten as
  \[\left\{ {{x}_{j}},{{p}_{k}} \right\}={{\delta }_{jk}}+{{x}_{j}}{{A}_{k}}+{{p}_{k}}{{b}_{j}}\]where $\left\{ \cdot ,\cdot \right\}$ is the GSPB.
\end{corollary}

\section{GCHS on the Riemannian manifold}
Let $M$ be a differentiable manifold of dimension $m$. A Riemannian metric on $M$ is a family of inner products $ g\colon T_{p}M\times T_{p}M\longrightarrow \mathbb{R},~~p\in M$ such that, for all differentiable vector fields $X,Y$ on $M$, $ p\mapsto g(X(p),Y(p))$ defines a smooth function $M\to \mathbb{R}$. the metric tensor can be written in terms of the dual basis $(dx_{1}, \cdots, dx_{n})$ of the cotangent bundle as ${ g=g_{ij}\,\mathrm {d} x_{i}\otimes \mathrm {d} x_{j}}$. Endowed with this metric, the differentiable manifold $(M, g)$ is a Riemannian manifold. In a local coordinate system $\left( U,{{x}_{i}} \right)$, connection $\nabla$ gives the Christoffel symbols, so now the structural derivative ${{A}_{i}}=\Gamma _{il}^{l}$ is now expressed as the special case of Christoffel symbols
\[
  \Gamma _{ij}^{k}=\frac{1}{2}{{g}^{kl}}\left( \frac{\partial {{g}_{jl}}}{\partial {{x}_{i}}}+\frac{\partial {{g}_{il}}}{\partial {{x}_{j}}}-\frac{\partial {{g}_{ij}}}{\partial {{x}_{l}}} \right)
\]
where $g^{ ij}$ is the inverse of $g_{ ij}$.
\begin{theorem}\label{th1}
  The GCHS on $(M, g)$ can be expressed
 \begin{align}
 \frac{\mathcal{D}{{x}_{k}}}{dt}
 & ={{J}_{kj}}\frac{\partial H\left( x \right)}{\partial {{x}_{j}}}+{{J}_{kj}}\Gamma _{ji}^{i}H\left( x \right)+{{x}_{k}}w\notag
\end{align}where $w={{J}_{ij}}\Gamma _{il}^{l}\frac{\partial H\left( x \right)}{\partial {{x}_{j}}}$ is the S-dynamics.
\begin{proof}
  Plugging the Levi-Civita connection into the GCHS \[\frac{\mathcal{D}{{x}_{k}}}{dt}={{J}_{kj}}{{\partial }_{j}}H+{{J}_{kj}}{{A}_{j}}H+{{x}_{k}}{{J}_{ij}}{{A}_{i}}{{D}_{j}}H\left( x \right)\]
then it proves the theorem.
\end{proof}
\end{theorem}
Clearly, theorem \ref{th1} has indicated GCHS as the real general dynamical system which has contained the new three parts only linked to the structural function, GHS is very ordinary part of GCHS.

\begin{corollary}\label{lem}
  S dynamics on $(M, g)$ is given by \[w=\hat{S}H\left( x \right)={{J}_{ij}}\Gamma _{li}^{l}\frac{\partial H\left( x \right)}{\partial {{x}_{j}}}\]where $\hat{S}={{J}_{ij}}\Gamma _{li}^{l}\frac{\partial }{\partial {{x}_{j}}}$ is structural operator.
\end{corollary}
Incidentally, from the corollary \ref{lem}, it is clear that S dynamics is related only to the structure function $s$, and TGHS is expressed in the following corollary
\begin{corollary}\label{lem1}
  TGHS on $(M, g)$ is given by
  \[\dot{x}_{k} =\frac{d{{x}_{k}}}{dt}={{J}_{kj}}{{\partial }_{j}}H+{{J}_{kj}}\Gamma _{lj}^{l}H\]
\end{corollary}
Clearly, corollary \ref{lem1} contains GHS as the rate of change of time, which is a whole. Only in this way can we describe the stability problems and development trends etc.
Therefore, only corollary \ref{lem1} can describe the Riemannian manifold with Levi-Civita connection, which can clearly describe the various properties of Hamiltonian dynamics on Riemannian manifold.

According to the corollary \ref{c3}, we can obtain a series of equations of GCHS on the Riemannian manifold

\begin{corollary}\label{c1}
The generalized force field $F_{k}$ on the Riemannian manifold is $$F_{k}={\dot{p}_{k}}=-{{D}_{k}}H$$where ${{D}_{k}}={{\partial }_{j}}+\Gamma _{lj}^{l}$.
\end{corollary}

\begin{corollary}
  By corollary \ref{c1}, canonical thorough Hamilton equations can be simply rewritten as
  \[{\dot{x}_{k}}={{J}_{jk}}{F_{j}},
        ~~{\dot{p}_{k}}=F_{k}\]on the Riemannian manifold.
\end{corollary}

\subsection{The geometrio on the Riemannian manifold}
According to theorem \ref{t2}, we give the geometrio in terms of the ${{x}_{k}},{{p}_{i}},H $ respectively on the Riemannian manifold
\begin{align}
 Geo\left( {{x}_{k}},{{p}_{i}},H \right) &=\hat{S}{{\left( {{x}_{k}},{{p}_{i}},H \right)}^{T}}=\left( \begin{matrix}
   {{J}_{ik}}\Gamma _{il}^{l}  \notag\\
   \Gamma _{il}^{l}  \\
   {{J}_{ij}}\Gamma _{il}^{l}{{\partial }_{j}}H  \notag\\
\end{matrix} \right) \\
 & ={{\left( {{b}_{k}},{{A}_{i}},w \right)}^{T}} \notag
\end{align}
Since ${{\left\{ s,f \right\}}_{GPB}}={{J}_{ij}}{{\partial }_{i}}s{{\partial }_{j}}f={{J}_{ij}}{{\partial }_{j}}f{{\partial }_{i}}s$, thusly, then the S-dynamics is
\begin{align}
 w & =\hat{S}H={{J}_{ij}}\Gamma _{il}^{l}{{\partial }_{j}}H={{\left\{ s,H \right\}}_{GPB}}=-{{\left\{ H,s \right\}}_{GPB}} \notag\\
 & =-{{J}_{ij}}{{\partial }_{i}}H\Gamma _{jl}^{l}={{J}_{ji}}{{\partial }_{i}}H\Gamma _{jl}^{l} \notag
\end{align}
Thus, then the geometrio in terms of the ${{x}_{k}},{{p}_{i}},H $ respectively on the Riemannian manifold is rewritten as
\begin{align}
 Geo\left( {{x}_{k}},{{p}_{i}},H \right) &=\hat{S}{{\left( {{x}_{k}},{{p}_{i}},H \right)}^{T}}={{\left( {{b}_{k}},{{A}_{i}},w \right)}^{T}} \notag\\
 & =\left( \begin{matrix}
   {{J}_{ik}}  \\
   1  \\
   {{J}_{ji}}{{\partial }_{i}}H  \notag\\
\end{matrix} \right)\Gamma _{il}^{l} \notag
\end{align}
As mentioned, it's clear to see that the geometrio on the Riemannian manifold is uniquely induced by the structure derivative or the special case of Christoffel symbols $\Gamma _{il}^{l}$. It evidently proves that the GCHS on the Riemannian manifold is certainly determined by the Christoffel symbols.

\subsection{Acceleration on Riemannian manifold}
\begin{definition}[Acceleration]\label{af}
The acceleration on generalized Poisson manifold $\left( P,S,\left\{ , \right\} \right)$ is defined as
\begin{equation}\label{eq22}
  a=\frac{{{\mathcal{D}}^{2}}x}{d{{t}^{2}}}=\ddot{x}
+2w{\dot{x}}+x\beta,~~~x\in {{\mathbb{R}}^{m}}
\end{equation}
where $$\ddot{x}
=\frac{{{d}^{2}}x}{d{{t}^{2}}},~~\beta ={{w}^{2}}+\frac{dw}{dt}$$ and ${\dot{x}}=-J\left( x \right)F$, Its component expression is  $${{a}_{i}}=\ddot{x}_{i}
+2w\dot{x}_{i}+{{x}_{i}}\beta $$

\end{definition}
Obviously, definition \ref{af} also reveals the inevitable problems and critical limitations of GHS. Conversely, GCHS as a universal and intact theoretical system can remedy the defects of GHS.  More importantly, acceleration on generalized Poisson manifold is directly and completely derived from GCHS, it has achieved the self consistency and original intention between the Hamiltonian structure and classical mechanics along with the GCHS, the theories are compatible with each other. This vital point goes far beyond GHS reaches.
\begin{theorem}\label{1a}
 Acceleration  on $(M, g)$ of GCHS  is
\[{{a}_{k}}=\frac{{{d}^{2}}{{x}_{k}}}{d{{t}^{2}}}+2w{{J}_{kj}}{{\partial }_{j}}H+2w{{J}_{kj}}\Gamma _{jl}^{l}H+{{x}_{k}}{{w}^{2}}+{{x}_{k}}\frac{dw}{dt} \]
where $\ddot{x}_{k}=\frac{d}{dt}\dot{x}_{k}$ is the second order derivative of time.
\end{theorem}Classically, there is only one term $\frac{d}{dt}\left( {{J}_{kj}}{{\partial }_{j}}H \right)$ for the GHS to depict the acceleration, it can't describe the Riemann manifold with Levi-Civita connection as mechanical system.

\section{Geospin matrix and S-dynamics}
In this section. Let us outline how one ascertains the invariance of the framework of the GCHS under the GSPB given in theorem \ref{t1}. Thus, this will provide some useful hints for the nature of the GCHS that is covariant and complete for the curved spacetime by putting the GCHS and geodesic equation together.

One wants to argue that the GCHS and geodesic equation are invariant under the coordinate transformation. This relies on the following one fact that is proven in a basic theoretic framework of the GCHS. More generally, we refer to any function $f$ satisfying $\left\{ f,H \right\}=0$ which exactly and similarly corresponds to the geodesic equation, techniques of the GCHS can be applied to examine the geodesic equation on Riemannian manifold.
\[\begin{matrix}
   \left\{ f,H \right\}=\frac{df}{dt}+wf=0  \\
   \frac{dv}{dt}+Wv=0  \\
\end{matrix}\]
Two equations have a nice geometric interpretation which can be deduced from the following simple analogy, the resulting differential equation are
\begin{equation}\label{eq12}
  \left\{ \begin{matrix}
   \frac{df}{dt}=-wf  \\
   \frac{dv}{dt}=-Wv  \\
\end{matrix} \right.
\end{equation}
Since $f$ is real-valued, $w$ is assumed to be real as well. Using each of this fact one can easily conclude that the GCHS is invariant under the coordinates transformation. Actually, the general solution to the ODE can be easily obtained. Note that the nature of the S-dynamic and the geospin matrix is the same as each other. More precisely, let $f=v$ be formally given intentionally:
\[\frac{dv}{dt}+wv=0,~~~ \frac{dv}{dt}+Wv=0\]
the former equation is built on the GCHS, meanwhile, the latter originates from the Riemannian geometry which is working for general relativity, namely, the geodesic equation. Thus, both equations hold the same mathematical expression. By comparison, it implies that the S-dynamics $w$ holds the same dynamical meanings as the geospin matrix $W$ does.

\subsection{Geospin matrix and S-dynamics on Riemannian geometry}
This section will analyze geospin matrix and S-dynamics on Riemannian geometry, we prove that geospin matrix and S-dynamics describe same thing for rotations.

To better realize this point, we can formally rewrites both dynamical system on Riemannian geometry together
\begin{align}
  & w=\hat{S}H ,~~ W=\left( W_{j}^{k} \right)\notag
\end{align}
where geospin variables $W_{j}^{k}=\Gamma _{ij}^{k}{{v}^{i}}$ and the corollary \ref{lem} on Riemannian geometry.  According to the S-dynamics expressed by \eqref{eq7}, that is, $w={\dot{x}_{i}}{{A}_{i}}$, As stated previously, the structural derivative on Riemannian geometry is the special case of Christoffel symbols given by
${{A}_{i}}=\Gamma _{il}^{l}$, then thusly, $w={\dot{x}_{i}}\Gamma _{il}^{l}$.
In particular, seen from \cite{3,5}, that ${{w}^{\left( j \right)}}=W_{j}^{j}={{A}_{i}}{\dot{x}^{i}}$ holds, if rewriting the S-dynamic is $w={\dot{x}^{i}}{{A}_{i}}$, namely, ${\dot{x}_{i}}$ is replaced by ${\dot{x}^{i}}$, then the formula of the S-dynamic $w$ equals the geospin variable ${{w}^{\left( j \right)}}=w$.

\subsection{The GCHS for geodesic equation}
In this section, we will use the GCHS to derive the equations like geodesic equation, and then we prove that the geodesic equation can be obtained based on the GCHS.

In the beginning, to consider the GCHS of the velocity fields ${{v}^{p}}=d{{x}^{p}}/dt,~~{{v}_{p}}=d{{x}_{p}}/dt$, it yields
\begin{align}\label{eq8}
 \frac{\mathcal{D}}{dt}{{v}^{p}} & =\left\{ {{v}^{p}},H \right\}={{J}_{ji}}{{F}_{j}}{{\partial }_{i}}{{v}^{p}}+{{v}^{p}}w \\
 & ={\dot{x}_{i}}\left( {{\partial }_{i}}{{v}^{p}}+{{v}^{p}}{{A}_{i}} \right) \notag={\dot{x}_{i}}{{D}_{i}}{{v}^{p}} \notag
\end{align}
where $w={{J}_{ji}}{{A}_{i}}{{F}_{j}}={\dot{x}_{i}}{{A}_{i}}$,  and $\frac{\mathcal{D}}{dt} ={\dot{x}_{i}}{{D}_{i}}$.
Subsequently, the GCHS in terms of the velocity field ${{v}^{p}}$ follows \[\frac{\mathcal{D}}{dt}{{v}^{p}}=d{{v}^{p}}/dt+w{{v}^{p}}\]
and
\begin{align}
\frac{\mathcal{D}}{dt}{{v}_{p}}  &=\left\{ {{v}_{p}},H \right\}={{J}_{ji}}{{F}_{j}}{{\partial }_{i}}{{v}_{p}}+{{v}_{p}}w \notag\\
 & ={\dot{x}_{i}}\left( {{\partial }_{i}}{{v}_{p}}+{{v}_{p}}{{A}_{i}} \right) \notag\\
 & ={\dot{x}_{i}}{{D}_{i}}{{v}_{p}}\notag
\end{align}
Similarly, the GCHS in terms of the velocity field ${{v}_{p}}$  is given by
\[\frac{\mathcal{D}}{dt}{{v}_{p}}=d{{v}_{p}}/dt+w{{v}_{p}}\]
Both equations can be united as a vector form
\[\frac{\mathcal{D}}{dt}v=dv/dt+wv\] by using $v=\left( {{v}^{p}} \right)=\left( {{v}_{p}} \right)$.
Actually, if ${{D}_{i}}$ is replaced by covariant derivative \eqref{eq11}
${{\nabla}_{i}}$, let's see how \eqref{eq8} becomes,
\[\frac{\mathcal{D}}{dt}{{v}^{p}}={\dot{x}_{i}}\left( {{\partial }_{i}}{{v}^{p}}+W_{i}^{p} \right)=d{{v}^{p}}/dt+W_{i}^{p}{{v}_{i}}\]
Then the geodesic equation reappears at $\frac{\mathcal{D}}{dt}{{v}^{p}}=0$, namely, like geodesic equation Eq\eqref{eq3},  $d{{v}^{p}}/dt+W_{i}^{p}{{v}_{i}}=0$.

As corollary \ref{c4} stated,  the covariant equilibrium equation of the GCHS in terms of the velocity vector is $\frac{\mathcal{D}}{dt}v=0$, that is, $dv/dt+wv=0$, clearly, $W\simeq w$, where $\simeq$ means equivalence effects.  In this way, we can say that the GCHS can naturally derive the geodesic equation, and GCHS self-consistently gives the explanation for the geodesic equation as general relativity stated.
\begin{theorem}
   The S-dynamic and the geospin variable satisfies $ w \simeq W$.
\end{theorem}
Obviously, the S-dynamic and the geospin variable compatibly remains the same completely on Riemannian geometry,  it clearly say that the covariant geodesic equation emerges in the framework of the GCHS, and   therefore, it well explains how both dynamical system describe one mode of motion.  In another way,  according to the GHS \eqref{eq4}, S-dynamics should be the form like $w={{v}_{i}}\Gamma _{il}^{l}$, now, we can evidently see that the formula form is the same on Riemannian geometry, it directly proves that the S-dynamic and the geospin matrix describe the same physical phenomena, meanwhile, the geodesic equation is covariant, thus, it reveals that the GCHS is also covariant as we name it in \cite{1}. We have obtained a compatible theory, and then it proves the GCHS along with the GSPB are consistent with the description of the non-Euclidean space, as a complete theory for general case.
\begin{corollary}
  The S-dynamic $w$ and the geospin variable $W$ on Riemannian geometry are respectively shown as
   $$w={\dot{x}_{i}}{{A}_{i}}\simeq W=\left( W_{j}^{k} \right)$$
   where ${{A}_{i}}=\Gamma _{il}^{l}$ is Christoffel symbols. Especially, ${{w}^{\left( j \right)}}=W_{j}^{j}\bowtie w$, where $\bowtie$ represents  the equivalence.
\end{corollary}
We can conclude that the scalar function $w$ for describing the S-dynamics in essence is same pattern as the geospin matrix $W$ created for the geodesic equation.  We know that the geodesic equation is induced and used for the curved spacetime, in particular, in general relativity, therefore, it also says that the frame structure of the GCHS is a real system for non-Euclidean case.    As a result, this indicates that the framework of the GCHS is a complete and compatible theory for general case.
Actually, it proves that we use different method to obtain compatible theory.
\begin{theorem}
  S-dynamics $w=\hat{S}H$ of the GCHS is the geometric angular frequency while geospin matrix $W$ for geodesic equation is for the geometric angular frequency, thus, the $W=\left( W_{j}^{k} \right)$ and $w=\hat{S}H$ are only for describing the curved spacetime.
\end{theorem}
As geodesic equation stated in general relativity, through comparative analysis, we can certainly conclude that the GCHS is self-contained and suitable for the general case. Especially, the GCHS can be naturally describing the general relativity.

\section{Conclusions}
Conclusively,  wa can say that the framework of the GCHS with the GSPB is complete and compatible with the curved spacetime, as stated, it naturally meets the non-Euclidean case. Undoubtedly, the GSPB as a mathematical tool can be used in other fields.  Then we prove that
\begin{align}
  Geo\left( {{x}_{k}},{{p}_{i}},H \right)&=\hat{S}{{\left( {{x}_{k}},{{p}_{i}},H \right)}^{T}}={{\left( {{b}_{k}},{{A}_{i}},w \right)}^{T}} \notag\\
 & ={{\left( {{J}_{ik}},1,{{J}_{ij}}{{\partial }_{j}}H \right)}^{T}}{{\partial }_{i}}s \notag
\end{align}Subsequently, we put it into the on Riemannian geometry.  By the way, we found that geospin matrix and S-dynamics describe the similar physical mode, exactly, it's for depicting some kinds of the rotation. It perfectly proves original intention of building the GCHS theory that is correct for a real world. Generally, the S-dynamics $w=\hat{S}H$ of the GCHS and geospin matrix $W=\left( W_{j}^{k} \right)$ are only used to the non-Euclidean case in essence.


\end{document}